\documentclass[a4paper]{amsart}
\usepackage{graphicx} 
 \usepackage{mathrsfs}
\usepackage{amsfonts}
\usepackage[active]{srcltx}
\usepackage{indentfirst,latexsym,bm,amsmath,pstricks,amssymb,amsthm,amscd}
\usepackage[all]{xy}

\numberwithin{equation}{section}

\def\cO{{\mathcal O}}

\newtheorem{theo}{Theorem}
\newtheorem{lem}[theo]{Lemma}

\newtheorem{coro}[theo]{Corollary}

\newtheorem{remk}[theo]{Remark}

\newtheorem{conj}[theo]{Conjecture}

\begin{document}
\title[A remark on a conjecture of Schnell]{A remark on a conjecture of Schnell }
 \author[Jun Lu]{Jun Lu}
 \address{Address of Jun Lu:~~~~School of Mathematical Sciences,  Key Laboratory of MEA(Ministry of Education) $\&$ Shanghai Key Laboratory of PMMP,  East China Normal University, Shanghai 200241, China}
\email{jlu@math.ecnu.edu.cn }

\author[Wan-Yuan Xu]{Wan-Yuan Xu*}
 \address{Address of Wan-Yuan Xu:~~~~ Department of Mathematics, Shanghai Normal University, Shanghai 200234, China}
 \email{wanyuanxu@shnu.edu.cn}

\footnotetext[1]{ \ {\itshape 2020 Mathematics Subject
Classification.} 14D06, 14J27, 14C20}

\footnotetext[2]{\ {\itshape Key words and phrases.} fibration, effective, pseudo-effective}

\footnotetext[3]{\ * Corresponding author}

 \begin{abstract}
 In this paper, we prove a conjecture of Schnell in the surface case.

 \end{abstract}
 
 \maketitle
 
\section{Introduction}
Let $f: X \to Y$ be a fibration between smooth projective varieties defined over an algebraically closed field $k$ with general fiber $F$. Here, by a fibration we mean $f$ is surjective with connected fibers such that $f_{*}\mathcal{O}_X=\mathcal{O}_Y$,
where $\mathcal{O}_X$ and $\mathcal{O}_Y$ represent the structure sheaves of $X$ and $Y$, respectively. 
We denote the canonical divisor of $X$ by $K_X$. The main purpose of this paper is to study a conjecture proposed by Schnell \cite[Conjectute 10.1]{CP11}.

\begin{conj}\label{SC}
   Let $f: X \to Y$  be a fibration as before with $\kappa(F)\geq 0$ and $H$ be an ample divisor on $Y$. Suppose that the divisor class $mK_X-f^*H$ is pseudo-effective for some $m\geq 1$, then $mK_X-f^*H$ is effective for $m$ sufficiently large and divisible.
    \end{conj}

In \cite{CP11}, Schnell proved that \textbf{Conjecture \ref{SC}}, together with the non-vanishing conjecture, is equivalent to the well known \textbf{Campana-Peternell Conjecture}. 

\begin{conj}[Campana-Peternell]\label{CP}
Suppose that for some $m\geq 1$, the divisor class $mK_X-f^*H$ is pseudo-effective, where $H$ is an ample divisor on $Y$. Then
$\kappa(X)\geq \dim Y$.
    \end{conj}
Hence, we believe that \textbf{Conjecture \ref{SC}} is an important part of solving the \textbf{Campana-Peternell Conjecture}.

\begin{remk}
    Here we formulate the \textbf{Campana-Peternell Conjecture} in an equivalent form as discussed in \cite[Section 4--6]{Schnell}. Although the original \textbf{Conjecture \ref{SC}} and the \textbf{Campana-Peternell Conjecture} are formulated for fibrations defined over the complex number field $\mathbb{C}$, one can ask the same question for fibrations defined over any algebraically closed field $k$. 
\end{remk}

Note that over $\mathbb{C}$, the \textbf{Campana-Peternell Conjecture} is true for $\dim X\leq 3$, since good minimal models exists up to dimension 3 \cite[p. 53--54]{CP11}. We also note that it is related to the behavior of Kodaira dimension in algebraic fiber spaces (see \cite{P21, PS18, PS22} and the references therein).

In \cite{Schnell}, by using recent deep theory about singular metrics on pluri-adjoint bundles, Schnell \cite[Theorem 12.1]{Schnell} proved \textbf{Conjecture \ref{SC}} over $\mathbb{C}$ under an additional assumption that the canonical divisor $K_Y$ of $Y$ is pseudo-effective. Furthermore, applying this result, he reduced the proof of the \textbf{ Campana-Peternell Conjecture} over $\mathbb{C}$ (modulo the non-vanishing conjecture) to the case where $Y$ is rationally connected. Unfortunately, he was not able to say anything about Conjecture \ref{SC} and \ref{CP} even in the case $Y=\mathbb{P}^1$ \cite[section 20]{Schnell}.

\begin{remk}
 It is easy to see that when $f$ is a surface fibration defined over $\mathbb{C}$, in order to prove conjecture \ref{SC} and \ref{CP} , it suffices to consider the case $Y=\mathbb{P}^1$. However, when $f$ is defined over an algebraically closed field $k$ with $\text{char} k=p>0$, we need to consider the base $Y$ is an arbitrarily smooth projective curve, since Schnell's proof for $g(Y)\geq 1$ is analytic. 
\end{remk}

In this paper, by using the surface fibration theory, we show that the \textbf{Conjecture \ref{SC}} is true for a surface fibration.
\begin{theo}\label{main}
   \textbf{Conjecture \ref{SC}} is true for a surface fibration defined over an algebraically closed field $k$. 
\end{theo}

\section{Proof of Theorem \ref{main}}

In this section, we prove our main Theorem \ref{main}.

 From now on, we assume that $f:X\to Y$ is a surface fibration of genus $g\geq 1$ defined over $k$ and $m_0K_X-f^*H$ is pseudo-effective for a given positive integer $m_0$. Let $b$ be the genus of $Y$ and let $H$ be an ample divisor on $Y$ with $l=\deg H$.

 For $r\gg 0$, $rH$ is very ample. Moreover, $(m_0r)K_X-f^*(rH)$ is also pseudo-effective and $$h^0(nK_X-f^*H)\geq h^0(nK_X-f^*(rH)) $$
 for $n\gg 0$. 
 Therefore, without loss of generality, we can assume that $H$ is a very ample and effective divisor in what follows. Thus one can write $f^*H=F_1+\dots+F_l$ where $F_i$'s are  general fibers of $f$.

Before proving the main theorem, we need some lemmas.

Let $D=\sum\limits_{i=1}^rn_i\Gamma_i$ be an effective divisor on $X$, where $n_i>0$ and $\Gamma_i$'s are the irreducible components of $D$. Let $\sigma: X\to X_0$ be the blowing-down contracting a $(-1)$-curve $E\subseteq X$ to a point $p\in X_0$. We assume that each $\Gamma_i$ cannot be contracted by $\sigma$ and define 
$D_0=\sum\limits_{i=1}^rn_i\overline{\Gamma}_i$, where each $\overline{\Gamma}_i$ is the image of $\Gamma_i$ under $\sigma$. 
One has 
$$\sigma^*D_0=D+\nu_0E,\quad \nu_0=\sum\limits_{i=1}^rn_i{\rm mult}_p(\overline{\Gamma}_i)\geq 0. $$
\begin{lem}\label{lem:H2H0-2}
 Let $D,D_0$ and $\sigma$ be as above. Then  $m_0K_X- D$ is pseudo-effective iff $m_0K_{X_0}- D_0$ is pseudo-effective. Moreover, $h^0(mK_X- D)=h^0(mK_{X_0} - D_0)$ for any $m>0$.   
\end{lem}
\begin{proof} 
It is easy to see that
$$mK_X-D=\sigma^*(mK_{X_0}-D_0)+(m+ \nu_0)E.$$
Recall that $\sigma$ is a contraction of a $(-1)$-curve $E$, it follows that $h^0(mK_X-D)=h^0(mK_{X_0} -D_0)$.

Suppose $(m_0K_X -D)$ is pseudo-effective, if $m_0K_{X_0} -D_0$ is not pseudo-effective, i.e., $(m_0K_{X_0}-D_0)H_0<0$ for some ample divisor $H_0\subseteq X_0$. We note that $\sigma^*H_0$ is nef, hence $(m_0K_X -D)\sigma^*H_0<0$, a contradiction. Theorefore, $m_0K_{X_0} -D_0$ is pseudo-effective.

Suppose $m_0K_{X_0} -D_0$ is pseudo-effective, take a very ample divisor $H$ in $X$(Without loss of generality, we can assume that $H\geq 0$). Let $H_0$ be the image of $H$ under $\sigma$. It is easy to see that $H_0$ is nef. So $(m_0K_X-D)H\geq (m_0K_{X_0}-D_0) H_0\geq 0$
\end{proof}

\begin{lem}\label{key-2}
    If  $m_0K_X-f^*H$ is pseudo-effective, then $\kappa(X)\geq 1$. Furthermore,  $h^0(nK_X-lF)>0$ for $n\gg 0$ whenever $\kappa(X)=2$.
\end{lem}
\begin{proof} 
 Let $\rho:X\to X_0$ be a sequence of blowing-downs and $X_0$ be the minimal model of $X$. Note that 
 $F_i$ cannot be contracted by $\rho$.  We take  $D_0=\sum\limits_{i=1}^l\overline{F}_i$, where each $\overline{F}_i$ is the image of $F_i$ under $\rho$.

By Lemma \ref{lem:H2H0-2}, $m_0K_{X_0}-D_0$ is pseudo-effective. So is $m_0K_{X_0}$. This implies that $\kappa(X)\ne -\infty$. 

    Suppose that $\kappa(X)=0$. It is well-known that $12K_{X_0}\equiv 0$ in this case \cite[Theorem 7.1]{L13}. 
    Thus  $12(m_0K_{X_0}-D_0)=-12D_0$ is not pseudo-effective, a contradiction.   Hence we have $\kappa(X)\geq 1$.

    If $\kappa(X)=2$, then 
    $$h^0(nK_X-f^*H) \geq h^0(nK_{X })-h^0(\cO_{f^*H}(nK_{X }))=bn^2-c,$$
    for some positive constants $b$ and $c$.
    Hnece, $h^0(nK_X-f^*H)>0$ for $n\gg 0$.
\end{proof}

\begin{coro}
Conjecture \ref{CP} is ture in the surface case.
\end{coro}
\begin{proof}
It is a direct corollary of Lemma \ref{key-2}.  
\end{proof}

\vspace{0.2cm}
From the above discussion, to prove Thoerem \ref{main}, it suffices to consider a fibration satisfying the following assumption.

{\bf Assumption:} {\itshape $f:X\to Y$ is a relatively minimal surface fibration such that  $\kappa(X)=1$ and $m_0K_X-f^*H$ is pseudo-effective.}
 
\begin{lem}
 If $h^0(mK_X-f^*H)>0$ for some $m>0$, then $h^0(nK_X-f^*H)>1$ for  $n\gg 0$. 
\end{lem}
\begin{proof}
It is obvious that 
 $$ h^0(nK_X-f^{*}H)\geq h^0((n-m)K_X)+h^0(mK_X-f^*H)-1>1,$$
 by $\kappa(X)=1$.
\end{proof}

Therefore, it suffices to find an integer $m>0$ such that $h^0(mK_X-f^*H)>0$.
 
Under our assumption, $X$ admits an elliptic fibration $h:X\to B$ (not necessary relatively minimal). Denote by $V$ the general fiber of $h$.
Since $m_0K_X-f^*H$ is pseudo-effective, one has
$$0\leq (m_0K_X-f^*H))V=-lFV\leq 0 $$ for a general fiber $F$ of $f$. 
So $FV=0$, namely, $h=f$. 

By our assumption, this elliptic fibration is relatively minimal. In this case, we have the canonical bundle formula for an elliptic fibration (cf. \cite[Theorem 2]{BM77} or \cite[Theorem 5.6]{L13})
$$K_X {\equiv} f^*D +\sum\limits_{i=1}^ka_iP_i,$$
where $V_i=m_iP_i (i=1,\cdots, k)$ are all the multiple fibers of $f$, $0\leq a_i<m_i, a_i\in \mathbb{Z}$ (``$\equiv$" means linear equivalent) and
$$\deg D=    \chi(\cO_X)+2b-2+{\rm length}(\mathcal{T}),$$
 where $\mathcal{T}$ is the torsion part of $R^1f_*\mathcal{O}_X$ (Note that $\mathcal{T}=0$ when $\text{char} k=0$).

\begin{lem}\label{lem:ell1-2}
    If 
    $$\deg D+\sum\limits_{i=1}^k\frac{a_i}{m_i}>0,$$
 then there exists an integer $m>0$ such that $h^0(mK_X-f^*H)>0$.   \end{lem}
\begin{proof}
Take $m=m_1m_2\cdots m_k  n$ for sufficiently divisible $n$.
Then 
$$mK_X\equiv f^*(mD)+ \sum\limits_{i=1}\frac{a_im}{m_i}V_i .$$
Note that for $n$ sufficiently large, one has
$$ m\deg D+\sum\limits_{i=1}\frac{a_im}{m_i}\gg l.$$
Let $p_i=f(V_i)$ and $q_j=f(F_j)$. 
Then we have  
$$h^0(mK_X-f^*H)=h^0(mD+\sum\limits_{i=1}^k\frac{a_im}{m_i}p_i-\sum\limits_{j=1}^lq_j)>0.$$
\end{proof}

\begin{proof}[{\bf Proof of Theorem \ref{main}}] 
Note that 
$$m_0K_X-f^*H\equiv_{\rm num} \left(m_0\deg D+\sum\limits_{i=1}^k\frac{a_im_0}{m_i} -l\right) F.$$ 
It is easy to see that
$m_0K_X-f^*H$  is pseudo-effective if and only if 
   \begin{equation}\label{eq1}
       \deg D+\sum\limits_{i=1}^k\frac{a_i}{m_i}\geq \frac{l}{m_0}(>0).
   \end{equation} 
    
Combining \eqref{eq1} with Lemma \ref{lem:ell1-2}, there exists an integer $m>0$ such that $h^0(mK_X-f^*H)>0$. This completes the proof.
\end{proof}

\noindent{\bf Acknowledgements}
This work is supported by NSF of China;  the first author is supported in part by Science and Technology Commission of Shanghai Municipality (No. 22DZ2229014). The second author is supported by Science and Technology Innovation Plan of Shanghai (Grant No. 23JC1403200).


\begin{thebibliography}{ABCD}

\bibitem{BM77} E. Bombieri, D. Mumford, {\it Enriques’ classification of surfaces in char. $p$, II}, in Complex analysis and algebraic geometry, Cambridge Univ. Press, 23--42, 1977.

 \bibitem{CP11} F. Campana and T. Peternell, {\it Geometric stability of the cotangent bundle and the universal cover of a projective manifold}, Bull. Soc. Math. France, 2011, {\bf{139} } no. 1, 41--74, with an appendix by Matei Toma.

\bibitem{L13} C. Liedtke, {\it  Algebraic surfaces in positive characteristic}, Birational geometry, rational curves, and arithmetic, New York, NY: Springer New York, 2013, 229--292.

\bibitem{P21} M. Popa, {\it Conjectures on the Kodaira dimension}, to appear in a volume in honor of V. Shokurov's 70th birthday, arXiv:2111.10900 (2021).

\bibitem{PS18} M. Popa and C. Schnell, {\it Viehweg’s hyperbolicity conjecture for families with maximal variation}, Invent. Math., (2017), {\bf208} no. 3,
677--713.

\bibitem{PS22} ------, {\it On the behavior of Kodaira dimension under smooth morphisms}, to appear in Algebraic Geometry,  arXiv:2202.02825v2.

\bibitem{Schnell} C. Schnell, {\it Singular metrics and a conjecture by Campana and Peternell},  preprint arXiv:2202.01295v1 (2022).


\end{thebibliography}
\end{document}